\documentclass[11pt,oneside,onecolumn,letterpaper]{article}
\usepackage{amsthm,amsmath,amssymb,euscript,amscd,subfigure,latexsym, psfrag, url}

\usepackage{epstopdf,graphicx}
\newtheorem{theorem}{Theorem}[section]
\newtheorem{corollary}[theorem]{Corollary}

\newtheorem{lemma}[theorem]{Lemma}

\newtheorem{remark}[theorem]{Remark}

\newcommand{\fdim}{r}

\begin{document}
\title{On the Damping-Induced Self-Recovery Phenomenon in \\ Mechanical Systems with \\ Several Unactuated Cyclic Variables}

\author{\\\\ \textbf{Dong Eui Chang} \\
        Department of Applied Mathematics\\ University of Waterloo\\
        Waterloo ON N2L 3G1 Canada,
        Tel: 519-888-4567 x37213, Fax: 519-746-4319 \\
        \url{dechang@uwaterloo.ca}\\\\
        \textbf{Soo Jeon}\\ Department of Mechanical and Mechatronics Engineering\\ University of Waterloo\\
        Waterloo ON N2L 3G1 Canada, Tel: 519-888-4567 x38898, Fax: 519-885-5862
        \\ \url{soojeon@uwaterloo.ca}\\
}
\date{}
\maketitle

\begin{abstract}
The damping-induced self-recovery phenomenon refers to the fundamental property of  underactuated mechanical systems: if an unactuated cyclic variable is under a viscous damping-like force and the system starts from rest,  then the cyclic variable  will always move back to its initial condition as the actuated variables come to stop. The regular momentum conservation phenomenon can be viewed as the limit of the  damping-induced self-recovery phenomenon in the sense that the self-recovery phenomenon  disappears  as the damping goes to zero. This paper generalizes the past result on damping-induced self-recovery for the case of a single unactuated cyclic variable to the case   of multiple unactuated cyclic variables. We characterize a class of external forces that induce  new conserved quantities, which we call the damping-induced momenta. The damping-induced momenta yield first-order asymptotically stable dynamics for the unactuated cyclic variables under some conditions, thereby inducing the self-recovery phenomenon. It is also shown that the viscous damping-like forces impose bounds on the range of trajectories of the unactuated cyclic variables. Two examples are presented to demonstrate the analytical discoveries: the planar pendulum with gimbal actuators and the three-link planar manipulator on a horizontal plane.

\end{abstract}

\section{Introduction}
Mechanical systems of kinematically coupled structures, or multibody systems, have a variety of engineering applications including robotic manipulators, manufacturing machines with articulated components, bipedal walking and musculoskeletal system modeling \cite{Krishna89}. An interesting aspect of such systems is how they  behave when certain degrees of freedom are left unactuated. In dealing with such underactuated mechanical systems, we often explore their fundamental properties. For instance, if the unactuated variables are cyclic (i.e. do not appear in the Lagrangian), then the momentum map associated with these variables will be conserved due to symmetry \cite{MaRa95}. One of the most celebrated examples of the symmetry and conservation of momentum is the falling cat problem \cite{Kane69,Montgomery93,Mo06}: a  cat always lands on its feet after released upside down from rest, by deforming its body while keeping a zero angular momentum. There are variants of the falling cat such as Elroy's beanie \cite{MaMoRa90} and springboard divers \cite{Fr79}. Another popular example is the system of  a human sitting on a rotating stool holding a wheel  \cite{ChJe12}.  Sitting on the rotating stool, one  spins the wheel by his hand while holding it horizontally. A reaction torque will be created and initiate the rotating motion of the stool in the opposite direction. As long as the wheel is rotating, the stool keeps rotating. After some time, if the person applies a braking force  halting the wheel spin, then the stool will also stop. All of these phenomena follow from angular momentum conservation, which holds in the absence of  external forces. If there is an external force, then the momentum is not conserved any longer in general and the consequent motion of the system will be different from that in the case of momentum conservation.

We are here mainly  interested in the case where the external forces are linear in the velocity, i.e. viscous damping-like forces. This was well studied in \cite{ChJe12} for the case of one-dimensional symmetry, i.e., $S^1$ or $\mathbb R^1$ symmetry. According to \cite{ChJe12}, in the stool-wheel system with viscous damping friction on the rotation axis of the stool, the stool does not keep rotating but meets a bound in its motion even if the person on the stool keeps spinning the wheel. The larger the angular velocity of the spin of the wheel is, the more angle the stool rotates by, but in the end the stool meets a bound in its motion, not being able to keep rotating, as long as the angular velocity of the wheel is bounded. This bound on the motion of the stool is called the damping-induced bound. Another phenomenon, which is perhaps more remarkable, in this stool-wheel system with the friction is the following. After some time, if the person on the stool stops the spinning of the wheel, then the stool does not stop. Instead, it asymptotically converges back  to its initial position, making the same  number of  net rotations in the  opposite direction. This phenomenon is called damping-induced self-recovery. From this example, one can see that the viscous damping-like force plays a role of a restoring force, but it is different from the spring force. To the knowledge of the authors, an example of damping-induced self-recovery and boundedness was first reported by Andy Ruina at a conference \cite{Ru10}, where he shows a video of his experiment and provides  an intuitive proof of the phenomena.
The damping-induced self-recovery phenomenon, without damping-induced boundedness,  is also explained in \cite{Gregg10} for the case where the damping coefficient is constant, but the proof therein is not complete.
Both phenomena of  damping-induced self-recovery and damping-induced boundedness were rigorously  and completely  proved for the first time in \cite{ChJe12}, where the damping coefficient is allowed to be a function of the cyclic variable and may  take some negative values. The results in \cite{ChJe12} as well as \cite{Ru10,Gregg10} are for the case of one-dimensional symmetry.


This paper generalizes the results in \cite{ChJe12} to the case of higher-dimensional Abelian symmetry. We consider mechanical systems with several unactuated cyclic variables that are subject to external forces   linear in the velocity. The systems are allowed to have some actuated variables that are subject to control forces. We find conditions on the linear forces under which the damping-induced  phenomena of self-recovery and boundedness  occur for the cyclic variables. The analytical approach taken in this paper for the proof of these two phenomena is different from that in \cite{ChJe12}. The total order of $\mathbb R$ was implicitly made use of in \cite{ChJe12}, but in this paper we employ a Lyapunov-like function to prove  the damping-induced self-recovery  and boundedness for the multiple cyclic variables. We first characterize a class of damping-like forces that induce new conserved quantities which we call in this paper  the damping-added momenta, and then construct a Lyapunov-like function using the damping-added momenta to prove the existence of the phenomena of damping-induced self-recovery and boundedness. To illustrate the main results, we consider two nontrivial examples of multibody systems that possess multiple cyclic variables: the planar pendulum with gimbal actuators and the three-link planar manipulator on a horizontal plane.

\section{Main Results}\label{sec:description}

In this paper, the norm $\| \cdot   \|$ denotes the Euclidean norm for vectors and the corresponding induced norm for matrices. For an invertible matrix $A = (a_{ij})$, the $(i,j)$-th entry of its inverse matrix is denoted by $a^{ij}$. For a symmetric matrix $A$ we denote its positive semi-definiteness by $A \succeq 0$. For two symmetric matrices $A$ and $B$ of the same size, $A \succeq B$ means $A-B \succeq 0$.

\subsection{Equations of Motion and Damping-Added Momenta}

Let $Q = Q_1 \times Q_2$ be an $n$-dimensional configuration space, where $Q_1 = \mathbb R^r$ and $Q_2$ is a smooth manifold of dimension $n-r$. If a unit circle $S^1$ appears as a factor of $Q_1$,   it is replaced by $\mathbb R$, which does not impose any restrictions on  the description of the dynamics. Let $q = ( x,  y) \in \mathbb R^\fdim \times  Q_2$ denote coordinates for $Q$, where
\[
 x = (x^\alpha)= (x^1, \ldots, x^\fdim), \quad  y = (y^a)= (y^{\fdim+1}, \ldots, y^n)
\]
and
\[
q = (q^i)  = (q^1, \ldots, q^\fdim; q^{\fdim +1}, \ldots, q^n) = (x^1, \ldots, x^\fdim; y^{\fdim +1}, \ldots, y^n).
\]
For notational convenience, the following three groups of indices are used in this paper:
\[
\underbrace{ \underbrace{1, \ldots, \fdim}_{ \alpha, \beta, \gamma, \ldots };  \underbrace{\fdim +1, \ldots, n}_{ a, b, c, \ldots}}_{i,j,k ...}.
\]

Consider a mechanical system with the Lagrangian
\begin{equation*}\label{def:Lag}
L(q,\dot q) = \frac{1}{2}m_{\alpha \beta}\dot x^\alpha \dot x^\beta + m_{\alpha a} \dot x^\alpha \dot y^a + \frac{1}{2}m_{ab}\dot y^a \dot y^b - V(q)
\end{equation*}
which is the kinetic minus potential energy of the system. Here we follow the Einstein summation convention.  It is assumed that the mass matrix
\[
m = \begin{pmatrix}
m_{\alpha \beta} & m_{\alpha b} \\
m_{a \beta} & m_{ab}
\end{pmatrix}
\]
is symmetric and positive definite.


 We  make the following four assumptions on the system:
\begin{itemize}
\item [A1)] The variables $x^\alpha$'s are cyclic \cite{MaRa95}, i.e.
 \begin{equation*}\label{e:cyclic variable}
 \frac{\partial L}{\partial x^\alpha}=0
 \end{equation*}
for all $\alpha = 1, \ldots, r$.

\item [A2)] Controls $u_{a}$'s are given in the directions of $y^{a}$'s.

\item [A3)] Each cyclic variable  $x^\alpha$ is under a general damping-like force (i.e. linear in the velocity) described as $-k_{\alpha \beta}(x)\dot x^\beta$, where each  coefficient $k_{\alpha \beta}(x)$ is a  continuously differentiable  function of $x$.

\item [A4)] The  coefficients $k_{\alpha\beta}$ satisfy
\begin{equation}\label{sym:cond:k}
k_{\alpha \beta} = k_{\beta\alpha}
\end{equation}
and
\begin{equation}\label{int:cond:k}
\frac{\partial k_{\alpha\beta}}{\partial x^\gamma} = \frac{\partial k_{\alpha\gamma}}{\partial x^\beta}.
\end{equation}

\end{itemize}
By A1) -- A3), the equations of motion of the system are written as
\begin{align}
\frac{d}{dt}\frac{\partial L}{\partial \dot x^\alpha}  &= - k_{\alpha \beta}\dot x^\beta, \quad \alpha = 1, \ldots, \fdim \label{e:equation of motion:x}\\
\frac{d}{dt}\frac{\partial L}{\partial \dot y^a} - \frac{\partial L}{\partial y^a} &= u_a, \quad a =  \fdim+1, \ldots, n. \label{e:equation of motion:y}
\end{align}
In individual coordinates, \eqref{e:equation of motion:x} and \eqref{e:equation of motion:y} can be expressed as
\begin{align}
&m_{\alpha \beta}\ddot x^\beta + m_{\alpha b} \ddot y^b + [ij,\alpha]\dot q^i \dot q^j = -k_{\alpha \beta}\dot x^\beta\label{EL:coord:1}\\
&m_{a\beta}\ddot x^\beta + m_{ab} \ddot y^b   + [ij,a]\dot q^i \dot q^j  + \frac{\partial V}{\partial y^a} = u_a,\label{EL:coord:2}
\end{align}
where $[ij,k]$ denotes the Christoffel symbol
\[
[ij,k]= \frac{1}{2} \left ( \frac{\partial m_{ik}}{\partial q^j} + \frac{\partial m_{jk}}{\partial q^i} - \frac{\partial m_{ij}}{\partial q^k}\right )
\]
for $i,j,k = 1, \ldots, n$.


If there were no external forces, i.e. $k_{\alpha \beta} =0$ for all $\alpha, \beta$, then the momenta $\frac{\partial L}{\partial \dot x^\alpha}$'s would be the first integrals of the system (\ref{e:equation of motion:x}), but  due to the external forces $F_\alpha = -k_{\alpha\beta} \dot x^\beta$'s,  $\frac{\partial L}{\partial \dot x^\alpha}$'s are not conserved any more.  We will here show that there are new conserved quantities  associated with \eqref{e:equation of motion:x} in place of $\frac{\partial L}{\partial \dot x^\alpha}$'s.

\begin{lemma}\label{lemma:fund}
Assumption A4) entails the following:
\begin{enumerate}
\item For each $\alpha$, there is a function $h_{\alpha} : \mathbb R^\fdim \rightarrow \mathbb R$ such that
\begin{equation}\label{prop:h}
\frac{\partial h_{\alpha}}{\partial x^\beta} = k_{\alpha\beta}.
\end{equation}
\item For each $h_{\alpha}$, there is a function $U : \mathbb R^\fdim \rightarrow \mathbb R$ such that
\begin{equation}\label{prop:V}
\frac{\partial U}{\partial x^\alpha} = h_{\alpha}.
\end{equation}Namely, the damping coefficient matrix ($k_{\alpha \beta}(x)$) is the second-order derivative matrix of $U$.
\end{enumerate}
\begin{proof}
\begin{enumerate}
\item By \eqref{int:cond:k} and Ponicar\'{e}'s lemma, there is a function $h_\alpha$ that satisfies \eqref{prop:h}.

\item By \eqref{sym:cond:k} and \eqref{prop:h}, we have ${\partial h_{\alpha}}/{\partial x^\beta} = k_{\alpha\beta} = k_{\beta\alpha} = {\partial h_{\beta}}/{\partial x^\alpha}$.
Hence, by Ponicar\'{e}'s lemma, there is a function $U : \mathbb R^\fdim \rightarrow \mathbb R$ such that \eqref{prop:V} holds.
\end{enumerate}
\end{proof}
\end{lemma}

Using Lemma \ref{lemma:fund}, we now show that the dynamics   \eqref{e:equation of motion:x} have $r$ first integrals.
\begin{theorem}
Let $h_\alpha$'s be functions that satisfy (\ref{prop:h}). Then, the $r$ functions
\begin{equation}\label{cons:quantity}
\frac{\partial L}{\partial \dot x^\alpha} + h_\alpha
\end{equation}
or equivalently
\begin{equation}\label{first:int}
m_{\alpha \beta} \dot x^\beta + m_{\alpha a} \dot y^a + h_\alpha
\end{equation}
for $\alpha =1, \ldots, \fdim,$ are the first integrals of the dynamics \eqref{e:equation of motion:x}.
\begin{proof}
Differentiating \eqref{cons:quantity} with respect to $t$ along the trajectory of the system gives
\begin{align*}
\frac{d}{dt} \left (\frac{\partial L}{\partial \dot x^\alpha} + h_\alpha \right ) = \frac{d}{dt}\frac{\partial L}{\partial \dot x^\alpha} + \frac{\partial h_{\alpha}}{\partial x^\beta} \dot x^\beta = -k_{\alpha\beta} \dot x^\beta + k_{\alpha\beta}\dot x^\beta =0
\end{align*}
due to \eqref{e:equation of motion:x} and \eqref{prop:h}. This completes the proof.
\end{proof}
\end{theorem}
The vector-valued map $(q, \dot q) \mapsto (\frac{\partial L}{\partial \dot x^\alpha} + h_\alpha)$  shall be called the damping-added momentum map.

\subsection{Damping-Induced Self-Recovery and Damping-Induced Boundedness}\label{sec:self-recovery and boundedness}
The first integrals in  \eqref{cons:quantity} depend on the initial condition $(x(0), y(0), \dot x(0), \dot y(0))$. Taking an arbitrary initial condition and letting $\mu_\alpha$ be the  initial value of the corresponding first integral,    $\frac{\partial L}{\partial \dot x^\alpha} + h_\alpha$, we have
\begin{equation}\label{eq:cons}
m_{\alpha \beta} \dot x^\beta + m_{\alpha a} \dot y^a + h_\alpha (x) = \mu_\alpha
\end{equation}
for all $t \in \mathbb R$.
By Lemma \ref{lemma:fund}, there is a function  $U$ on $\mathbb R^\fdim$ such that \eqref{prop:V} holds. By adding $-\mu_\alpha x^\alpha$ to it, we can obtain a function $U_\mu : \mathbb R^\fdim \rightarrow \mathbb R$ such that
\begin{equation}\label{dVmu}
 \frac{\partial U_\mu(x)}{\partial x^\alpha} = h_\alpha(x) - \mu_\alpha.
\end{equation}
The conservation equations \eqref{eq:cons} for the damping-added momenta  can be regarded as   first-order differential equations for the  cyclic variables $x^\alpha$'s. The damping-induced self-recovery phenomenon is a direct consequence of asymptotic stability of these first-order dynamics under some conditions on the function $U_\mu$. Let us make the following assumptions on $U_\mu$ which constitute sufficient conditions for the self-recovery phenomenon to occur:
 \begin{itemize}
\item [A5)] The function $U_\mu$ has a unique critical point, denoted $x_e$, and it is a minimum point of $U_\mu$.

\item [A6)] There is a number $\delta_1 >0$ such that
\[
\inf_{\| x -x_e\|\geq \delta_1} U_\mu(x)  > U_\mu(x_e).
\]
\item [A7)] There is a number $\delta_2 >0$ such that
\[
\inf_{\| x -x_e\|\geq \delta_2} \| dU_\mu(x)\| >0,
\]
where $dU_\mu = (\frac{\partial U_\mu}{\partial x^1}, \ldots, \frac{\partial U_\mu}{\partial x^\fdim})$.
\end{itemize}
Assumption A6) guarantees that there does not exist any sequence $\{x_k\}$ in $\mathbb R^r$ with $\lim_{k\rightarrow\infty} \|x_k\| = \infty$ such that $\lim_{k\rightarrow \infty} U_\mu(x_k) = U_\mu(x_e)$.  Likewise, assumption A7) guarantees that there does not exist any sequence $\{x_k\}$ in $\mathbb R^r$ with $\lim_{k\rightarrow\infty} \|x_k\| = \infty$ such that $\lim_{k\rightarrow \infty} \|dU_\mu(x_k)\| = 0$.

  We now state one of the two  main theorems of this paper.
\begin{theorem}[Damping-Induced Self-Recovery]\label{thm:main:1}
Suppose that controls $u_a(t)$'s are chosen such that $\lim_{t\rightarrow \infty} \dot y(t) = 0$  and there are  numbers $c_1 >0$, $c_2>0$ and $c_3>0$ such that
\begin{align}
 &c_1I \preceq  (m_{\alpha \beta} (y(t)))  \preceq c_2I \label{c1:m:c2}\\
 &\| (m_{\alpha a} (y(t))) \| \leq c_3\label{m:c3}
\end{align}
for all $t\geq 0$.
Then, $\lim_{t\rightarrow \infty} x(t) = x_e$ and $\lim_{t\rightarrow \infty} \dot x(t) = 0$. In particular, if the initial condition is such that $\dot x(0) =0$ and $\dot y(0)=0$, then  $\lim_{t\rightarrow \infty} x(t) = x(0)$.
\begin{proof}
By \eqref{prop:h}, integration of \eqref{e:equation of motion:x} with respect to $t$ yields \eqref{eq:cons} for some constants  $\mu_\alpha$'s, which can be written as
\begin{equation}\label{dot:x:alpha:eq}
\dot x^\alpha =  - m^{\alpha\beta} \frac{\partial U_\mu}{\partial x^\beta} - m^{\alpha \beta} m_{\beta a} \dot y^a
\end{equation}
by \eqref{dVmu}.
By adding a constant, we may assume that $U_\mu (x_e) =0$.  For each $s>0$, define an open set $W_s$  in $\mathbb R^r$ by
\[
W_s = \{ x\in \mathbb R^\fdim \mid U_\mu (x) < s\}.
 \]
 Take any $\epsilon >0$. We will find a number $T>0$ such that $\| x(t) - x_e \| < \epsilon$ for all $t\geq T$.

 By A5) and A6), there is a sufficiently small $\delta >0$ such that
 \begin{equation}\label{Uimplyepsilon}
x \in W_\delta \Rightarrow \| x - x_e \| < \epsilon.
\end{equation}
Replacing $\delta$ by a smaller positive number if necessary,  we can find, by A5) -- A7), an $\epsilon_1 >0$ such that
 \begin{equation}\label{UdV}
x \notin W_\delta \Rightarrow \| dU_\mu (x) \| \geq \epsilon_1.
\end{equation}
Since the matrix $(m_{\alpha\beta})$ is symmetric and positive definite,  \eqref{c1:m:c2} implies
\begin{equation}\label{c2:c1:inverse}
\frac{1}{c_2}I \preceq  (m^{\alpha \beta} (y(t)))   \preceq \frac{1}{c_1}I.
\end{equation}
Choose a number $0<\ell <1$ such that
\begin{equation}\label{ineq:4:d}
\frac{1}{c_2} \left (1-\frac{1}{2}\ell \right )^2 - \frac{1}{4c_1}\ell^2 > 0,
\end{equation}
which is always possible since the left hand side of \eqref{ineq:4:d} is continuous in $\ell$ and is positive at $\ell=0$.
Since $\lim_{t\rightarrow \infty} \dot y(t) =0$ and the matrix $(m_{\alpha a}(y(t)))$ is bounded, there is a $T_1 >0$ such that
\begin{equation}\label{depsilon1}
\| (m_{\alpha a} (y(t)) \dot y^a(t)) \| <  \ell \epsilon_1
\end{equation}
for all $t \geq T_1$.  Whenever $x(t) \notin W_\delta$ for some $t \geq T_1$,  we have, by \eqref{dot:x:alpha:eq}, \eqref{UdV}, \eqref{c2:c1:inverse} and \eqref{depsilon1},
\begin{align}
\frac{dU_\mu(x(t))}{dt} &= \frac{\partial U_\mu(x)}{\partial x^\alpha}\dot x^\alpha  \nonumber \\
&= \frac{\partial U_\mu}{\partial x^\alpha} \left ( - m^{\alpha\beta} \frac{\partial U_\mu}{\partial x^\beta} - m^{\alpha \beta} m_{\beta a} \dot y^a\right ) \nonumber \\
&= -  m^{\alpha\beta} \left ( \frac{\partial U_\mu}{\partial x^\alpha} + \frac{1}{2}m_{\alpha a} \dot y^a \right )\left (  \frac{\partial U_\mu}{\partial x^\beta} + \frac{1}{2}m_{\beta b} \dot y^b \right ) + \frac{1}{4}m^{\alpha \beta} m_{\alpha a}\dot y^a m_{\beta b} \dot y^b \nonumber  \\
&\leq -\frac{1}{c_2} \left \| \left (  \frac{\partial U_\mu}{\partial x^\alpha} + \frac{1}{2}m_{\alpha a} \dot y^a \right )  \right \|^2 + \frac{1}{4c_1}  \| ( m_{\alpha a}\dot y^a) \|^2 \nonumber \\
&\leq -\frac{1}{c_2} \left |  \left \| dU_\mu\right \|  - \frac{1}{2} \| (m_{\alpha a}\dot y^a)\| \right |^2 + \frac{1}{4c_1} \ell^2\epsilon_1^2 \nonumber \\
&\leq -\frac{1}{c_2} \left | \epsilon_1 - \frac{1}{2}\ell\epsilon_1 \right |^2 + \frac{1}{4c_1} \ell^2\epsilon_1^2 \nonumber  \\
&= - \left ( \frac{1}{c_2} \left (1-\frac{1}{2}\ell \right )^2 - \frac{1}{4c_1}\ell^2 \right ) \epsilon_1^2, \label{dVmu7}
\end{align}
where the right-hand side of \eqref{dVmu7} is negative by \eqref{ineq:4:d}. Hence, $U_\mu(t)$ decreases at least linearly in time as long as $x(t) \notin W_\delta$. By definition of $W_\delta$,   there must exist $T > T_1$ such that $x(t) \in W_\delta $ for all $t\geq T$. Thus,   $\|x(t) - x_e\| < \epsilon$ for all $t \geq T$ by \eqref{Uimplyepsilon}. Therefore, $\lim_{t\rightarrow \infty } x(t) = x_e$. By taking the limit of both  sides of \eqref{dot:x:alpha:eq}, we obtain $\lim_{t\rightarrow\infty} \dot x(t) = 0$ since $dU_\mu (x_e) = 0$ by A5). In particular, if $\dot x(0) =0$ and $\dot y(0) = 0$, then $x_e = x(0)$ by A5) and \eqref{dot:x:alpha:eq}, so $\lim_{t\rightarrow \infty} x(t) = x(0)$.
\end{proof}
\end{theorem}

\begin{remark}The result in Theorem \ref{thm:main:1} is global. In order to get a local result, one has only to assume A5) -- A7) in a neighborhood of $x_e$ and to assume that the trajectory $x(t)$ stays in the neighborhood.
\end{remark}

As was discovered in the case with a single cyclic variable \cite{ChJe12}, the viscous damping force not only induces self-recovery but also imposes a bound to the range of the cyclic variable. Such a boundedness property also holds for multiple cyclic variables  as stated in the following theorem.
\begin{theorem}[Damping-Induced Boundedness]\label{thm:damping:bound}
Suppose that the function $U_\mu$ satisfies
\begin{equation}\label{dVmuinfty}
\lim_{\|x\| \rightarrow \infty} \| dU_\mu (x) \| = \infty
\end{equation}
and
\begin{equation}\label{U:infty}
\lim_{\|x\| \rightarrow \infty } U_\mu(x) = \infty.
\end{equation}
If controls $u_a(t)$'s are chosen such that $\dot y(t)$ is bounded and there exist $c_1 >0$, $c_2 >0$ and $c_3 >0$ such that (\ref{c1:m:c2}) and (\ref{m:c3}) hold for all $t\geq 0$, then $x(t)$ is bounded.
\begin{proof}
Since $\dot y(t)$ is bounded by assumption, there is a number $c_4>0$ such that $\| \dot y(t)\| \leq c_4$ for all $t\geq 0$.  By \eqref{dVmuinfty}, there is a number $\ell>0$  such that
\[
\|x - x_e\| \geq \ell \Rightarrow \| dU_\mu (x) \| \geq 1 + \frac{c_2c_3c_4}{c_1}.
\]
Suppose that $x(t)$ is not bounded. Then, there are   numbers $0<t_1 < t_2$ such that $\| x(t) - x_e\| \geq \ell$ for all $t\in [t_1, t_2]$ and $U_\mu(x(t_1)) < U_\mu(x(t_2))$ by (\ref{U:infty}). Then, for all $t \in [t_1, t_2]$
\begin{align*}
\frac{dU_\mu(t)}{dt}  &= \frac{\partial U_\mu}{\partial x^\alpha} \left ( - m^{\alpha\beta} \frac{\partial U_\mu}{\partial x^\beta} - m^{\alpha \beta} m_{\beta a} \dot y^a\right )\\
&\leq -\frac{1}{c_2} \| dU_\mu (x(t))\|^2 + \frac{c_3c_4}{c_1}\|dU_\mu (x(t))\| \\
&=  \left (-\frac{1}{c_2} \| dU_\mu (x(t))\| + \frac{c_3c_4}{c_1}  \right ) \|dU_\mu (x(t))\| \\
&\leq \left ( -\frac{1}{c_2} \left(1 + \frac{c_2c_3c_4}{c_1}\right) + \frac{c_3c_4}{c_1} \right ) \|dU_\mu (x(t))\|\\
&\leq -\frac{1}{c_2}.
\end{align*}
Hence,
\begin{align*}
0 < U_\mu(x(t_2)) - U_\mu(x(t_1)) = \int_{t_1}^{t_2} \frac{dU_\mu}{dt}(t) dt \leq -\frac{1}{c_2} (t_2 - t_1) < 0,
\end{align*}
which is a contradiction. Therefore, $x(t)$ is bounded.
\end{proof}
\end{theorem}

\subsection{Diagonal Damping Force: A Special Case}
Suppose that there are continuous functions $k_1, \ldots, k_\fdim : \mathbb R \rightarrow \mathbb R$ such that   the damping coefficients $k_{\alpha \beta}$ are given as
\begin{align}
&k_{11} (x) = k_1(x^1); \quad k_{22}(x) = k_2(x^2); \quad \cdots \quad; \quad k_{\fdim\fdim}(x) = k_\fdim (x^\fdim); \label{k:alphas}\\
&k_{\alpha\beta}(x) = 0 \quad \textup{ for $\alpha \neq \beta$}. \nonumber
\end{align}
Notice that we do not assume that $k_\alpha(x^\alpha)$'s  take only non-negative values though we call them damping coefficients for convenience.


A function $h_\alpha$ satisfying \eqref{prop:h} is given by
\begin{equation}\label{simple:h}
h_\alpha(x^\alpha) = \int_0^{x^\alpha} k_\alpha(s) ds,
\end{equation}
where there is no summation over the index $\alpha$. Given $\mu = (\mu_1, \ldots, \mu_r) \in \mathbb R^r$, the function $U_\mu$  defined by
\begin{equation}\label{simple:Vmu}
U_\mu (x) = \sum_{\alpha=1}^\fdim \int_0^{x^\alpha} (h_\alpha(s) - \mu_\alpha) ds
\end{equation}
satisfies \eqref{dVmu}.

\begin{corollary}[Damping-Induced Self-Recovery]\label{cor:self}
Suppose that the functions $k_\alpha$'s given in (\ref{k:alphas}) and the functions $h_\alpha$'s defined in \eqref{simple:h} satisfy the following:
\begin{itemize}
\item [(i)] For each $\alpha$, the equation $h_\alpha (s) - \mu_\alpha = 0$ has a unique root, which is denoted by $x_e^\alpha$.

\item [(ii)] For each $\alpha$, $k_\alpha (x^\alpha_e) >0$.

\item [(iii)] For each $\alpha$, there is an open interval $I$ containing $x_e^\alpha$ such that
\[
\inf_{s \in \mathbb R\backslash I } |h_\alpha(s) - \mu_\alpha | >0.
\]
\end{itemize}
Then, the function $U_\mu$ defined in \eqref{simple:Vmu} satisfies A5) -- A7) such that the conclusions in Theorem \ref{thm:main:1} hold true.

\begin{proof}
Take any $\alpha$ between $1$ and $\fdim$. Since $h^\prime_\alpha(x_e^\alpha) = k_\alpha(x_e^\alpha) >0$, the function $h_\alpha (s) - \mu_\alpha$ is increasing over an open interval containing $x_e^\alpha$. Thus, by condition (i), we have $h_\alpha (s) - \mu_\alpha >0$ for all $s > x_e^\alpha$ and $h_\alpha (s) - \mu_\alpha <0$ for all $s < x_e^\alpha$. It is now straightforward to show that the function $U_\mu$ defined  in \eqref{simple:Vmu} satisfies A5) -- A7).
\end{proof}
\end{corollary}

\begin{corollary}[Damping-Induced Boundedness]
Suppose that the functions $h_\alpha$ defined in \eqref{simple:h} satisfy the three conidtions in Corollary \ref{cor:self} and the following condition:
\begin{itemize}
\item [(iv)] For each $\alpha$, $\lim_{s\rightarrow \infty} h_\alpha(s) = \infty$ and $\lim_{s\rightarrow -\infty} h_\alpha(s) =- \infty$.
\end{itemize}
Then, the conclusion in Theorem \ref{thm:damping:bound} holds true.
\begin{proof}
It is easy to show that $U_\mu$ given in \eqref{simple:Vmu} satisfies \eqref{dVmuinfty} and \eqref{U:infty}.
\end{proof}
\end{corollary}

\subsection{Choice of Control for  Damping-Induced Self-Recovery}
In Theorems \ref{thm:main:1} and \ref{thm:damping:bound} we have assumed that an appropriate control law exists to satisfy some conditions on the trajectory $y(t)$. We now constructively show that such a control law exists.

\begin{lemma}[\cite{Spong94}]\label{lemma:spong}
The  control law
\begin{equation}\label{control:ua}
u_a = f_a(y,\dot x, \dot y) + g_{ab}(y) \tau^b,
\end{equation}
where $\tau^b$'s are new control variables and
\begin{align*}
f_{a} &= \left( [ij,a] - m_{a\alpha}[ij,\beta]m^{\alpha\beta} \right )\dot q^i \dot q^j -  m_{a\alpha} m^{\alpha\beta}k_{\beta \gamma} \dot x^\gamma+ \frac{\partial V}{\partial y^a},\\
g_{ab}&= m_{ab} -m_{a\alpha}m_{b\beta} m^{\alpha\beta},
\end{align*}
transforms the system in (\ref{EL:coord:1}) and (\ref{EL:coord:2}) to the following system:
\begin{align}
\dot x^\alpha &=  - m^{\alpha\beta}\frac{\partial U_\mu}{\partial x^\beta} - m^{\alpha\beta} m_{\beta b} \dot y^b\label{seqn:1}\\
\ddot y^a &= \tau^a, \label{seqn:2}
\end{align}
where $U_\mu$ is a function that satisfies \eqref{dVmu}.
\begin{proof}
 Solve  \eqref{EL:coord:1} for $\ddot x^\beta$,  substitute it into  \eqref{EL:coord:2}, and apply the control in  \eqref{control:ua} to obtain  \eqref{seqn:2}. Since $dh_\alpha/dt = k_{\alpha \beta} \dot x^\beta$, integration of \eqref{e:equation of motion:x} with respect to $t$ yields \eqref{eq:cons} where $\mu_\alpha$ is the value of the first integral.
Solving \eqref{eq:cons} for $\dot x^\alpha$ gives \eqref{seqn:1}.
\end{proof}
\end{lemma}

Equation \eqref{seqn:2} implies that we have full control over the motion of $y^a$ variables. Hence, it is always possible to find a control law to satisfy the assumptions in Theorems \ref{thm:main:1} and \ref{thm:damping:bound}; refer to   \cite{ChJe12} for more details on the method of choosing controls $u_a$'s. For example,  suppose that  $y_{\rm d}(t) = (y^a_{\rm d}(t))$ is a reference trajectory that $y(t)$ must follow. Then, we can choose the following control law for $\tau^a$:
\begin{equation}\label{e:error dynamics from PD}
\tau^a = {\ddot {y} }_{\rm d}^a(t) - c^a_1 ( \dot y^a - {\dot  {y}}_{\rm d}^a (t) ) - c^a_0 (  y^a-  {y_{\rm d} ^a} (t))
\end{equation}
with $c^a_1 >0$ and $c^a_0 >0$  such that the tracking error $e^a(t):= y^a(t) - y^a_{\rm d}(t)$ obeys the following exponentially stable dynamics:
\[
\ddot e^a  + c^a_1 \dot e^a + c^a_0 e^a  =0.
\]
Thus, $y^a(t)$ and $\dot y^a(t)$ converge exponentially to the reference trajectory $y^a_{\rm d}(t)$ and $\dot y^a_{\rm d}(t)$, respectively,  for each $a = r+1, \ldots, n$.

\section{Examples}\label{sec:examples}

In this section, we take two examples of mechanical systems with multiple cyclic variables to demonstrate the phenomena of damping-induced self-recovery and boundedness.

\subsection{Planar Pendulum with Gimbal Actuators}
\begin{figure}[!htp]
      \centering
      \includegraphics[trim = 0mm 0mm 0mm 0mm, clip,scale=0.42]{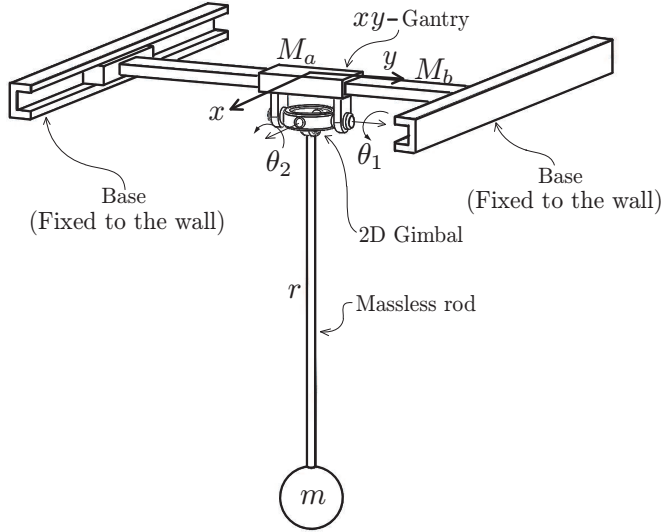}
      \caption{The planar pendulum with gimbal actuators.}
      \label{f:planar pendulum}
\end{figure}
Consider the planar pendulum system in Fig. \ref{f:planar pendulum}. The planar motion of the base block is unactuated and constrained to move freely along $x$ and $y$ directions only. The pendulum rod is assumed to be actuated by the gimbal-like mechanism which can apply torques along $\theta_1$ and $\theta_2$ rotational axes. Choose  coordinates as  $(q^1, q^2; q^3, q^4) = (x,y; \theta_1, \theta_2)$, where $x$ and $y$ are unactuated but subject to viscous damping-like forces and $\theta_1$ and $\theta_2$ are actuated by controls (i.e. $n=4$ and $r=2$). The mass matrix ${M}$ is given by
\begin{equation*}
{M}=\begin{pmatrix}M_a+M_b+m & 0 & -mrc^{\phantom{1}}_1c^{\phantom{1}}_2 & mrs^{\phantom{1}}_1s^{\phantom{1}}_2\\
0 & M_a+m & 0 & mrc^{\phantom{1}}_2\\
-mrc^{\phantom{1}}_1c^{\phantom{1}}_2 & 0 & mr^2c_2^2 & 0\\mrs^{\phantom{1}}_1s^{\phantom{1}}_2 & mrc^{\phantom{1}}_2 & 0 & mr^2\end{pmatrix},
\end{equation*}where $c_i$ and $s_i$ denote the $\cos{\theta_i}$ and $\sin{\theta_i}$, respectively and the parameters in ${M}$ are listed in Table  \ref{t:planar pendulum}. Since $x$ and $y$ do not show up in ${M}$ and also are not affected by the gravity, they are cyclic variables.
\begin{table}[!htp]
\begin{center}\caption{Parameters for planar pendulum} \label{t:planar pendulum}
\begin{tabular}{p{1cm}p{3.7cm}cl}
\hline
$\phantom{s}$ & Parameter & Value & Unit\\
\hline
$M_a$ & Slider mass & 2 & [kg]\\
$M_b$ & Gantry bar mass & 3 & [kg]\\
$m$ & Pendulum ball mass & 3 & [kg]\\
$r$ & Rod length  & 0.5 & [m] \\
\hline
\end{tabular}
\end{center}
\end{table}

\begin{figure}[!htp]
      \centering
      \includegraphics[trim = 0mm 0mm 0mm 0mm, clip,scale=0.75]{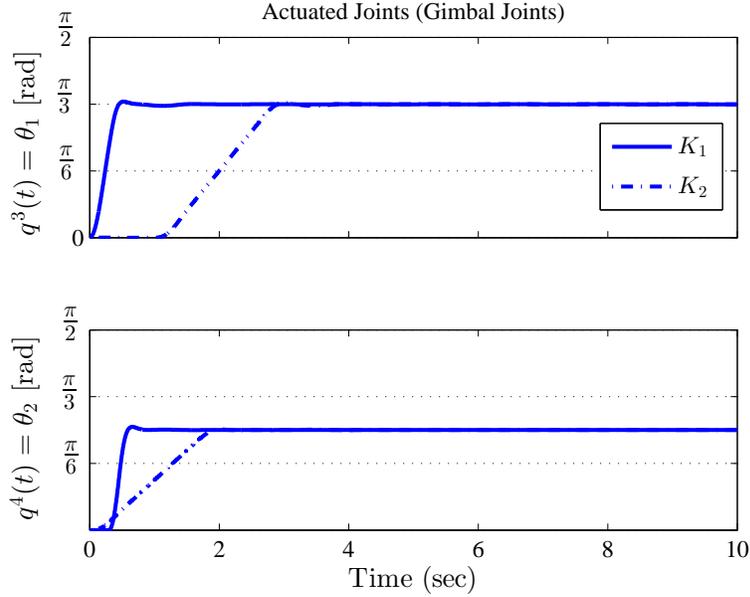}
      \caption{Time trajectories of the actuated variables, $\theta_1(t)$ and $\theta_2(t)$.}
      \label{f:planar_pendu_th}
\end{figure}

To demonstrate the self-recovery phenomenon in this case, we apply the control law in Lemma \ref{lemma:spong} to  control $\theta_1$ and $\theta_2$ such that they move from 0 [rad] to $\frac{\pi}{3}$ [rad] and $\frac{\pi}{4}$ [rad], respectively, within 3 seconds; see Fig. \ref{f:planar_pendu_th}. For the unactuated cyclic joints, $x$ and $y$, we simulate two different cases of damping matrix $K = (k_{\alpha\beta})_{1\leq \alpha, \beta \leq 2}$:
\[
K_1 = \begin{pmatrix}3 & 0\\0 & 3\end{pmatrix},\quad K_2 = \begin{pmatrix}5 + 2\cos{x} & 4\\4 & 4+ 2\cos y\end{pmatrix}.
\]
Notice that  $K_1$ and $K_2$ are the second-order derivative matrices of $U_1(x,y) = \frac{3}{2}\left(x^2 + y^2\right)$ and $U_2(x,y) = \frac{1}{2}x^2 + 2(x+y)^2  - 2\cos x -2\cos y$, respectively. One can easily verify that both $U_1$ and $U_2$ satisfy assumptions A5) -- A7) and equations (\ref{dVmuinfty}) and (\ref{U:infty}).

\begin{figure}[!htp]
\centering
\subfigure[Time trajectories.]
{\label{f:planar_pendu_xy_time}
\includegraphics[trim = 0mm 0mm 0mm 0mm, clip,scale=0.7]{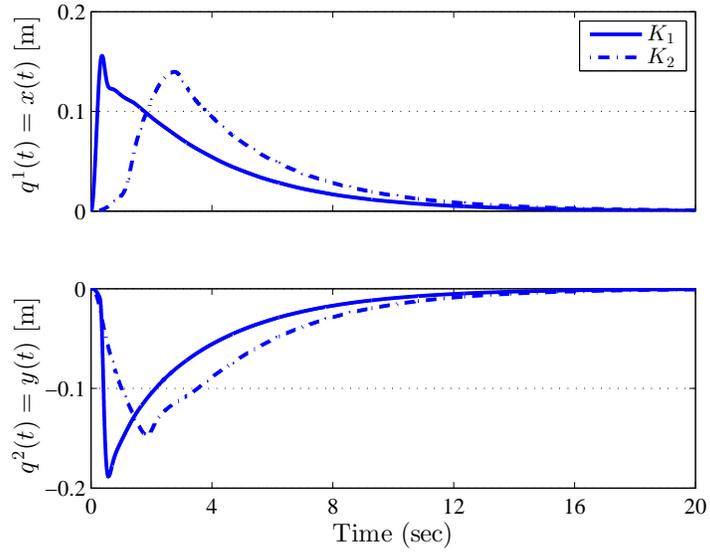}}
\subfigure[$x-y$ plot.]
{\label{f:planar_pendu_xy}
\includegraphics[trim = 0mm 0mm 0mm 0mm, clip,scale=0.7]{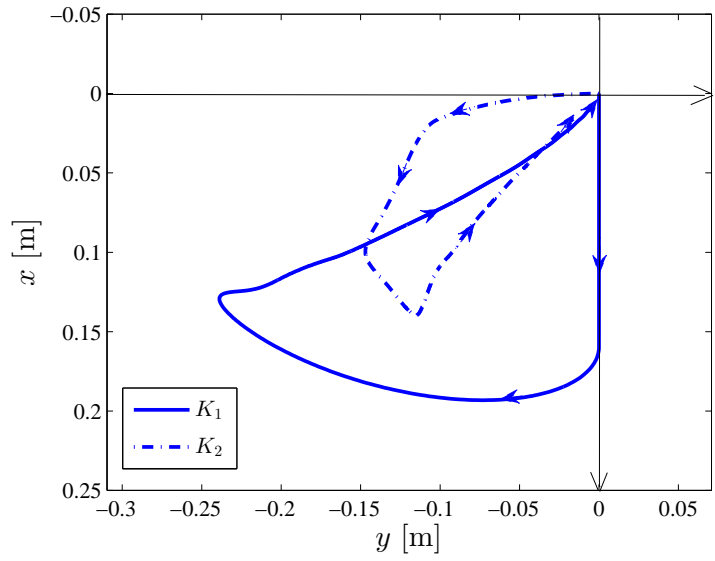}}
\caption{Trajectories of unactuated cycle variables, $x(t)$ and $y(t)$.}
\label{f:planar_pendu_xy_and_time}
\end{figure}


The self-recovery phenomenon for both cases is verified in Fig. \ref{f:planar_pendu_xy_and_time}. From Fig. \ref{f:planar_pendu_xy_time} and Fig. \ref{f:planar_pendu_th}, one can see that $x(t)$ and $y(t)$ converge to the origin, which is the initial condition,  after $\theta_1 (t)$ and $\theta_2 (t)$ settle down. The planar motion is described in Fig. \ref{f:planar_pendu_xy}, where we can clearly see the pendulum base automatically returns to its initial position for both cases. 



\subsection{Three-Link Manipulator on a Horizontal Plane}
\begin{figure}[!htp]
      \centering
      \includegraphics[trim = 0mm 0mm 0mm 0mm, clip,scale=0.42]{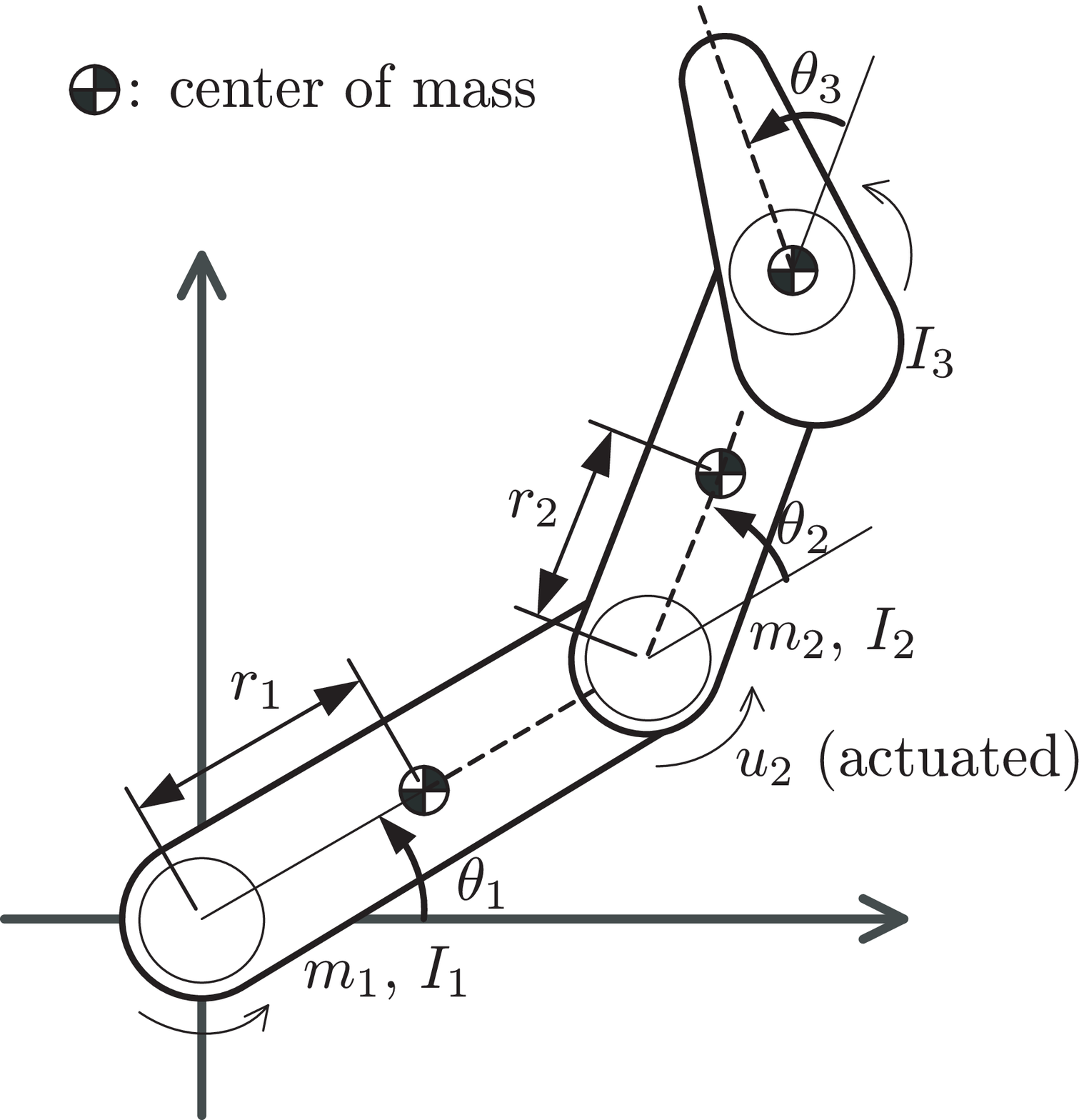}
      \caption{A three-link manipulator on a horizontal plane with $\theta_1$ and $\theta_3$ unactuated.}
      \label{f:threelink}
\end{figure}


Consider the three-link open chain manipulator moving on a horizontal plane (i.e. no gravity effect) in Fig.  \ref{f:threelink}. One can easily see that the first joint angle $\theta_1$ is a cyclic variable, which is always the case for planar kinematic chains whose motion is constrained to a horizontal plane. For this particular system in  Fig.  \ref{f:threelink},  the third joint angle $\theta_3$ is also cyclic since  the center of mass of the third link is located on  its axis of rotation. Thus, we anticipate the self-recovery effect will occur in these two joint variables. To be compatible with the notation introduced in Section \ref{sec:description}, we choose the configuration variables as $q^1 = \theta_1$, $q^2 = \theta_3$ and $q^3 = \theta_2$. Then, the mass matrix is written as
\begin{equation*}
{M}=\begin{pmatrix}\alpha + 2\beta c_2 & I_3 & \delta + \beta c_2\\ I_3 & I_3 & I_3 \\ \delta + \beta c_2 & I_3 & \delta \end{pmatrix},
\end{equation*}where $c_2$ denotes $\cos{\theta_2}$, and parameters $\alpha$, $\beta$ and $\delta$ are given by
\begin{equation*}
\begin{split}
\alpha &= I_1 + I_2 + I_3 + m_1r_1^2 + m_2(\ell_1^2+r_2^2) + m_3(\ell_1^2+\ell_2^2),\\
\beta &= \ell_1(m_2r_2 + m_3\ell_2),\\
\delta &= I_2 + I_3 + m_2r_2^2 + m_3\ell_2^2,
\end{split}
\end{equation*}
where individual parameter values are listed in Table \ref{t:three link}.  We use the following damping coefficient matrix:
\[
K = \begin{pmatrix}6 & 0\\0 & 3\end{pmatrix}
\]
which is the second-order derivative matrix of the function $U(\theta_1,\theta_3) = \frac{3}{2}\left(2\theta_1^2 + \theta_3^2\right)$.

\begin{table}[!htp]
\begin{center}\caption{Parameters for the three-link manipulator} \label{t:three link}
\begin{tabular}{p{1cm}p{5.5cm}cl}
\hline
$\phantom{s}$ & Parameter & Value & Unit\\
\hline
$\ell_1$ & Link 1 length  & 0.5 & [m] \\
$\ell_2$ & Link 2 length  & 0.5 & [m] \\
$r_1$ & Location of link 1 center of mass  & 0.1 & [m] \\
$r_2$ & Location of link 2 center of mass  & 0.1 & [m] \\
$I_1$ & Link 1 moment of inertia & 2 & [$\textrm{kg}\cdot\textrm{m}^2$]\\
$I_2$ & Link 2 moment of inertia & 2 & [$\textrm{kg}\cdot\textrm{m}^2$]\\
$I_3$ & Link 3 moment of inertia & 2 & [$\textrm{kg}\cdot\textrm{m}^2$]\\
$m_1$ & Link 1 mass & 10 & [kg]\\
$m_2$ & Link 2 mass & 10 & [kg]\\
$m_3$ & Link 3 mass & 10 & [kg]\\
\hline
\end{tabular}
\end{center}
\end{table}

\begin{figure}[!htp]
      \centering
      \includegraphics[trim = 0mm 0mm 0mm 0mm, clip,scale=0.75]{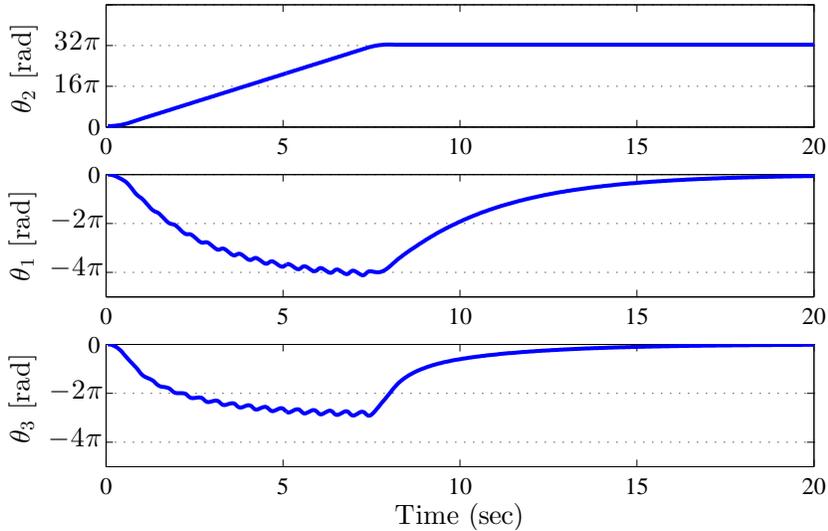}
      \caption{Time trajectories of joint angles, $\theta_2(t)$ (actuated), $\theta_1(t)$ (unactuated),  and $\theta_3(t)$ (unactuated).}
      \label{f:threelink_th}
\end{figure}

Figure \ref{f:threelink_th} plots the time trajectories of the joint angles. The top plot shows the controlled motion of the second joint $\theta_2$ while the second and the third  present the time trajectories of the two unactuated joints $\theta_1$ and $\theta_3$, respectively. A tracking controller is designed according to \eqref{e:error dynamics from PD} such that the actuated second joint $\theta_2$  makes sixteen full revolutions (i.e. $32\pi$ [rad]) during the first 8 seconds of time. The initial values for the joint angles are set to $\theta_2(0)=\frac{\pi}{3}$, and $\theta_1(0)=\theta_3(0)=0$. Figure \ref{f:threelink_th} clearly shows that both of the unactuated cyclic variables $\theta_1$ and $\theta_3$ converge back to their initial positions right after $\theta_2$ is regulated at $32\pi$ [rad]. Note that $\theta_1$ returns back to its original position after having made two full turns, demonstrating that the self-recovery effect is global. Besides the self-recovery phenomenon, the trajectories of $\theta_1$ and $\theta_3$ in Fig. \ref{f:threelink_th} also demonstrate the damping-induced boundedness derived in Theorem \ref{thm:damping:bound}. More specifically, we can observe that $\theta_1$ and $\theta_2$ oscillate around $-4\pi$ and $-3\pi$, respectively while $\theta_2$ is still increasing.

%

\section{Conclusion and Future Work}\label{sec:conclusions}
The phenomena of damping-induced self-recovery and damping-induced boundedness have been generalized to mechanical systems with several unactuated cyclic variables. The major contribution of this paper comes from characterizing the  damping coefficient matrix as the second-order derivative matrix of a function and identifying a class of such functions that  guarantee the self-recovery and boundedness of all  unactuated cyclic variables. Regular momentum conservation is the limit of the damping-induced self-recovery  as the damping disappears, in the sense that the recovery phenomenon vanishes in the limit. Non-trivial examples of mechanical systems with multiple cyclic variables are provided to demonstrate the theoretical discoveries.

It will be interesting to study the possibility of occurrence of damping-induced self-recovery and boundedness for the case of non-Abelian symmetry in the unactuated variables. We are currently examining the system of a spacecraft with internal rotors, where the spacecraft experiences viscous damping friction as it rotates.

\end{document}